\documentclass{article}
\usepackage{latexsym,amssymb}
\usepackage{amsfonts}
\usepackage{amsthm}
\newtheorem{theorem}{Theorem}[section]
\newtheorem{lemma}[theorem]{Lemma}
\theoremstyle{proposition}
\theoremstyle{definition}

\newtheorem{example}[theorem]{Example}
\newtheorem{proposition}[theorem]{Proposition}

\theoremstyle{remark}

\begin{document}
\title{Distances Induced by Barbilian's Metrization
Procedure}
\author{Wladimir G. Boskoff, Marian G. Ciuc\u{a},
Bogdan D. Suceav\u{a}}
\date{}
\maketitle

\begin{abstract}
Several authors have pointed out the connection between
Barbilian's metric introduced in 1934 and the recent study of
Apollonian metrics. We provide examples of various distances that
can be obtained by Barbilian's metrization procedure and we
discuss the relation between this metrization procedure and
important Riemannian and generalized Lagrangian metrics (in the
sense presented in \cite{A1996, MAB}). Then we prove an extension
of Barbilian's metrization procedure.  \end{abstract}

\section{Introduction: The Metrization
Procedure}

Barbilian's metrization procedure was introduced in \cite{B1934}
and it was the subject of an inspiring correspondence between D.
Barbilian and W. Blaschke \cite{B1940} in 1934 and thereafter. The
theory received  a larger audience due to P. Kelly \cite{K1954}
and a major development due to D. Barbilian
\cite{B1959a,B1959b,B1960,BR1962}. Over the years, the paper
\cite{B1934} has been cited many times. Recent studies are due to
A.F. Beardon \cite{Beardon1998},  F. Gehring and K. Hag
\cite{GH2000}, as well as P. H\"ast\"o, Z. Ibragimov and other
authors \cite{H2003, H2003b, H2004b, H2004, H2005b,
H2005,HL2004,HPS2005, I2002, I2003a, I2003b}.
 The geometric viewpoint is discussed in the
monograph \cite{WB1996carte}. All of these works cite and have a
common source in Barbilian's paper \cite{B1934}. The examples
explored in the present work aim to discuss Barbilian's
metrization procedure in the context of its relations with various
classes of metrics, as for example Riemann, Finsler, Lagrange or
Lagrange generalized metrics (see \cite{BCS2000, MAB}).

The following construction is given by Barbilian \cite{B1959a} and
it is the development of the idea from \cite{B1934}. Consider two
arbitrary sets $K$ and $J.$ The function $f : K \times J
\rightarrow \mathbb{R}_+^*$ is called an influence  of the set $K$
over $J$ if for any $A,B \in J$ the ratio
$g_{AB}(P)=\frac{f(P,A)}{f(P,B)}$ has a maximum $M_{AB} \in
\mathbb{R}$ when $P \in K.$ Note that $g_{AB}:K \rightarrow
\mathbb{R}_+^*.$ In \cite{B1959a} it is pointed out that if we
assume the existence of $\max g_{AB}(P),$ when $P \in K,$ then
there also exists $m_{AB} = \min_{P \in K} g_{AB}(P) = \frac{1}{
M_{BA}}.$

For example (see \cite{B1959a}), if $T$ is a topological space,
$K$ a compact subset in $T$, and $J$ some arbitrary subset,then
any function $f: K \times J \rightarrow \mathbb{R}_+^*$ continuous
in the first argument is an influence on $J.$ It is known since
\cite{B1959a} that $d: J \times J \rightarrow \mathbb{R}_+$ given
by\begin{equation}\label{metrica1}d(A,B) = \ln \frac{\max_{P\in K}
g_{AB}(P)}{\min_{P\in K} g_{AB}(P)}\end{equation}is a
semidistance, i.e.: (1) if $A=B$ then $d(A,B)=0;$ (2) $d$ is
symmetric; (3) $d$ satisfies triangle inequality.

The influence $f: K \times J \rightarrow \mathbb{R}_+^*$ is called
effective if there is no pair $(A,B) \in J \times J$ such that the
ratio $g_{AB}(P)=\frac{f(P,A)}{f(P,B)}$ is constant for all $P \in
K.$ In \cite{B1959a} it is shown that if $f: K \times J
\rightarrow \mathbb{R}_+^*$ is an effective influence, then
(\ref{metrica1}) is a distance.

\section{Examples}
\begin{example} Barbilian's
metrization procedure yields the  Euclidean distance in a plane $(
\pi )$ in $\mathbb{R}^3,$ if we consider a plane $(\delta)$
parallel to the plane $(\pi)$ and take $J=(\pi), K=(\delta)$, and
the influence function
 $f: K \times J \rightarrow \mathbb{R}_+^*$, $f(M,A)=
\exp \circ \left[ \frac{1}{2} || (Pr \times Id) (M,A)|| \right] =
e^{\frac{1}{2}||M'A||}.$\end{example}

\begin{example}
Barbilian's metrization procedure yields the spherical distance in
a complete sphere in $\mathbb{R}^3.$ \end{example}

To see this, consider two concentric spheres $S_1$ and $S_2$ in
$\mathbb{R}^3,$ and let their common center be $O.$ We take
$S_1=K$ and $S_2=J,$ and  $A,B \in J$ and $M\in K.$ Denote by
$\{M'\}=(OM \cap J$ and define $Pr$ the radial projection from
$S_1$ to $S_2$ given by $Pr(M)=M'.$ Denote by $( \ . \ )$ the
spherical distance, and consider the influence function $f: K
\times J \rightarrow \mathbb{R}_+^*$, $f(M,A)=\exp \circ
\left[\frac{1}{2} ( (Pr \times Id) (M,A) ) \right] =
e^{\frac{1}{2}(M'A)}.$

Thus, Barbilian's metrization procedure can generate Riemannian
metrics. Our goal is to show that Barbilian's metrization
procedure generates, for other choices of $K,J,$ and $f$, Lagrange
generalized metrics not reducible to a Riemannian, Finslerian or
Lagrangian metric.

To complete our discussion, we mention here the following result,
needed in the remaining part of this section. This is a particular
form of the result from \cite{B1960}, part 2, paragraph 7, and a
version of the argument used in \cite{H2004} in the proof of Lemma
3.5.
\begin{lemma}\label{mimetic} Let $K$ and $J$ be two
subsets of the Euclidean plane $\mathbb{R}^2,$ and $K = \partial
J.$ Consider the influence $f(M,A)= ||MA||,$ where by $||MA||$ we
denote the Euclidean distance. Consider $$g_{AB} (M) =
\frac{f(M,A)}{f(M,B)} = \frac{||MA||}{||MB||}$$ and consider the
distance induced on $J$ by the Barbilian's metrization procedure,
$d^B(A,B).$ Suppose furthermore that for $M \in K$ the extrema
$\max g_{AB}(M)$ and $\min g_{AB}(M)$ for any $A$ and $B$ in $J$
are attained each in an unique point in $K.$ Then:

(a) For any $A \in J$ and any line $d$ passing through $A$ there
exist exactly two circles tangent to $K$ and to $d$ in $A.$

(b) The metric induced by the Barbilian distance has the form
\begin{equation}\label{ecuatia_mimetica}ds^2 = \frac{1}{4}
\left(\frac{1}{R}+\frac{1}{r} \right)^2 (dx_1^2 +
dx_2^2),\end{equation} where $R$ and $r$ are the radii of the
circles described in (a).\end{lemma}

\begin{example} Let $K=$ be the line $\{
y=0 \}$ in the $xy$-plane. Let $J= \{ (x,y) / y>0 \}.$ Take the
function $||MA||$ as influence. Then the associated ratio is
$f(M)=\frac{||MA||}{||MB||}.$ By applying Barbilian's metrization
procedure, we only need to analyze the existence of minimum and
maximum for the function$$g(x) = \frac{x^2 - 2x_0 \cdot x + x_0^2
+ y_0^2}{x^2 - 2 x_1 \cdot x + x_1^2 +y_1^2}.$$ A straightforward
application of Lemma \ref{mimetic} yields, after computations $$R
= \frac{y \sqrt{m^2+1}}{-1+\sqrt{m^2+1}}$$ and $$r =
\frac{y\sqrt{m^2+1}}{1+\sqrt{m^2+1}},$$ that is $\frac{1}{4}
\left(\frac{1}{R}+\frac{1}{r} \right)^2=\frac{1}{y^2},$ i.e.
$ds^2=\frac{1}{y^2} (dx^2+dy^2),$ which is the Poincar\'e metric
on the upper half-plane.\end{example}

\begin{example}Consider
$\mathbb{R}^2$ endowed with the Euclidean distance $||.||$. It is
known from \cite{B1959a} that for any circle $K$ of radius $\rho$
in $\mathbb{R}^2$, and for $J$ the interior of $K$, a Barbilian's
distance is obtained in $J$ by taking the influence $f(P,A)=
||PA||.$ For a given point $(x,y)$ in J and for an arbitrary line
of slope $m$ passing through $(x,y),$ we find $$R =
\frac{\sqrt{m^2+1}}{2} \cdot \frac{\rho^2 - x^2 -
y^2}{\rho\sqrt{m^2+1} - x m +y}.$$ Similarly, we get, $$r =
\frac{\sqrt{m^2+1}}{2} \cdot \frac{\rho^2 - x^2 -
y^2}{\rho\sqrt{m^2+1} + x m -y}.$$ Hence, we proved the metric
relation$$\frac{1}{4} \left( \frac{1}{R} + \frac{1}{r} \right)^2 =
\frac{4 \rho^2}{(\rho^2-x^2 - y^2)^2}.$$ By a straightforward
computation, we can easily see that the Gaussian curvature of this
metric is $\kappa _g=-1.$ Therefore this Riemannian metric
generates the hyperbolic geometry on the disk.
\end{example}

For the next example, we apply Lemma \ref{mimetic} to the
following.

\begin{proposition}Barbilian's metrization
procedure on $$K=\{(x,0) \in \mathbb{R}^2 / x>0 \} \cup \{(0,y)\in
\mathbb{R}^2 /y>0 \},$$$$J= \{(x,y)\in \mathbb{R}^2 / x>0, y>0
\}$$for the influence $f: K \times J \rightarrow \mathbb{R}_+^*,$
given by $f(M,A) = ||MA||$, yields the metric that at $(x_0,y_0)
\in J$ satisfies\begin{equation}\label{metrica_neriemanniana}ds^2
=\frac{(y_0m+x_0+(x_0+y_0) \sqrt{m^2+1} )^2}{4x_0^2y_0^2
(m^2+1)}(dx^2 + dy^2),\end{equation}
$m=\frac{\dot{y}}{\dot{x}}|_{(x_0,y_0) },$ where the metric
(\ref{metrica_neriemanniana}) is a generalized Lagrange metric
that is not reducible to a Riemannian, Finslerian or Lagrangian
metric.
\end{proposition}

{\it Proof:} Denote as above $g_{AB}:K \rightarrow
\mathbb{R}_+^*,$ given by $$g_{AB}(M) =\frac{f(M,A)}{f(M,B)}=
\frac{||MA||}{||MB||}.$$First, we need to show that $g_{AB}$
admits maximum and minimum. Consider the points $A,B \in J,$ $M
\in K,$ and denote by $A_1$ the foot of perpendicular from $A$ to
the $y-$axis and by $A_2$ the foot of perpendicular from $A$ on
the $x-$axis. Consider the inversion centered in $A$ and of power
$||AA_1||^2.$ This inversion induces the correspondences $A_1
\rightarrow A_1,$ $A_2 \rightarrow A_2'$ (such that $A_2' \in
OA_2$ and$||AA_2|| \cdot ||AA_2'|| = ||AA_1||^2),$ $O \rightarrow
O'$ (such that $O' \in AO$ and $A_1O' \bot AO$), $B \rightarrow
B'$ such that $B' \in AB$ and $||AB|| \cdot ||AB'||= ||AA_1||^2$).
The positive part of the $y-$ axis is transformed in the arc of
circle $\mathcal{C}_1$ of endpoints $A$ and $O',$ and it is part
of the circle of diameter $AA';$ more precisely is the arc that
contains the point $A_1.$ The positive part of the $x-$ axis is
transformed in the arc of circle $\mathcal{C}_2$ of endpoints $A$
and $O',$ and it is part of the circle of diameter $AA_2'$, more
precisely the arc that contains the point $A_2.$ The inverse of a
point $M \in K$ is part of the union of the two arcs described
above. Keeping in mind that
\begin{equation}\label{formula}||B'M'|| = ||AA_1||^2 \cdot
\frac{||BM||}{||AM|| \cdot ||AB||}= \frac{||AA_1||^2}{||AB||}\cdot
\frac{||BM||}{||AM||},\end{equation} we get that $||B'M'||$ is
maximum whenever $\frac{||AM||}{||BM||}$ is minimum. Denote by
$M_1'$ the point on $\mathcal{C}_1 \cup \mathcal{C}_2$ for which
is attained the maximum of the Euclidean distance $||B'M'||.$ The
ray $AM_1'$ intersects $K$ in $M_1$ for which $$m
=\frac{||AM_1||}{||BM_1||} = \min_{M \in K}
\frac{||AM||}{||BM||}.$$ From (\ref{formula}) we deduce also that
there exists a point $M_2'$ for which $||B'M_2'||$ is the minimum
for $||B'M'||,$ when $M' \in \mathcal{C}_1 \cup \mathcal{C}_2.$
The inverse of $M_2'$ is $M_2$, obtained at the intersection
between $AM_2'$ and $K$ and it has the property
$$\mathcal{M}=\frac{||AM_2||}{||BM_2||}= \max_{M \in
K}\frac{||AM||}{||BM||}.$$ This allows us to conclude that the
formula $d^B(A,B)= \ln \frac{\mathcal{M}}{m}$ produces a Barbilian
distance in $J.$ Now we obtain the coefficients of the metric from
Lemma \ref{mimetic}. Consider the arbitrary point $A(x_0,y_0) \in
J$ and the line $(d)$ of equation $y-y_0=m(x-x_0).$ By Lemma
\ref{mimetic} there exist the circles $\Gamma_1$ and $\Gamma_2$
tangent to the line $d$ in $A$ and tangent to $K.$ Denote by
$O_1(x_1,y_1)$ the center of the circle $\Gamma_1$ and by
$O_2(x_2,y_2)$ the center of the circle $\Gamma_2.$ To determine
the rays of the two circles described in Lemma \ref{mimetic} (a)
we have the conditions
\begin{equation}
\label{prima}y_1 -y_0 = -\frac{1}{m}(x_1 - x_0),  \ \ \ x_1^2 =
(x_1-x_0)^2 + (y_1-y_0)^2,\end{equation} with $x_0 > x_1,$ and
\begin{equation}
\label{a doua} y_2-y_0= -\frac{1}{m}(x_2 - x_0), \ \ \ y_2^2 =
(x_2-x_0)^2 + (y_2 - y_0)^2,\end{equation} for $y_0 >y_1.$ From
(\ref{prima}) and (\ref{a doua}), respectively, we obtain:
\begin{equation}R_1=x_1=\frac{x_0
\sqrt{m^2+1}}{m+\sqrt{m^2+1}}, \ \ \ R_2=y_2=\frac{y_0
\sqrt{m^2+1}}{1+\sqrt{m^2+1}}.\end{equation} Therefore, by
applying Lemma \ref{mimetic} the metric is expressed as in
(\ref{metrica_neriemanniana}). For the directions
$m=\frac{\dot{y}}{\dot{x}}$ with $\dot{x} >0,$ the metric has the
coefficients
\begin{equation}g_{11}=g_{22}=
\frac{\dot{x}(y \cdot \dot{y}+ x \cdot \dot{x} +(x+y)
\sqrt{\dot{x}^2+\dot{y}^2})^2}{4xy(\dot{x}^2+\dot{y}^2)}, \ \ \
g_{12}=g_{21}=0.
\end{equation}
This metric (see \cite{M1994, MAB}) is a generalized Lagrange
metric, since the tensor expressed above is a $d$-tensor. To see
this, remark that the metric is $0-$homogeneous, and $\det g =
(g_{11})^2,$ therefore it is positive definite. According to
section 2.2 from \cite{MAB}, the metric
\ref{metrica_neriemanniana} is reducible to a Lagrangian metric if
and only if the Cartan tensor $C_{ijk}=\frac{1}{2}\frac{\partial
g_{ij}}{\partial x^k}$ is totally symmetric (see \cite{MAB},
section 4.1, Theorem 1.1.). The condition of symmetry reduces for
the metric (\ref{metrica_neriemanniana}) to $$\frac{\partial
g_{11}}{\partial \dot{y}}=\frac{\partial g_{12}}{\partial
\dot{x}}.$$However, $\frac{\partial g_{12}}{\partial
\dot{x}}\equiv 0$ and $\frac{\partial g_{11}}{\partial \dot{y}}
\neq 0,$ which proves that the Cartan tensor is not totally
symmetric. Therefore, the metric (\ref{metrica_neriemanniana}) is
not reducible to a Lagrangian metric. If the metric is not
reducible to a Lagrangian metric, it is not reducible to either a
Finslerian metric or a Riemannian metric. \qed

\section{An Extension of Barbilian's Metrization
Procedure}

Now we present an extension of Barbilian's metrization procedure.
Our motivation to produce this extension is the fact that in the
case when $K$ is a circle in the plane and $J$ is its interior, if
we remove one point $L$ from $K$, we can not apply the classical
Barbilian's metrization procedure considering the influence of $K
- \{ L \}$ over $J.$ Suppose that $K$ and $J$ are arbitrary sets
and that they satisfy the {\it general extremum requirement}, that
is for any $A$ and $B$ in $J$ it exists  $\sup g_{AB}(Q) <
\infty$, when $Q \in K.$ As we have seen in the case of maximum,
if there exists $\sup_{P\in K} g_{AB}(P) < \infty$ then there
exists $\inf_{P\in K} g_{AB}(P)$ and it equals $[\sup_{P\in K}
g_{BA}(P)]^{-1}.$  We have the following (see also \cite{WB2002},
p.10).

\begin{theorem} Suppose that $g$
satisfies the general extremum requirement. Then the function
$d^s: J \times J \rightarrow \mathbb{R}_+$ given by$$d^s(A,B) =
\ln \frac{\sup_{P\in K} g_{AB}(P)}{\inf_{P\in K} g_{AB}(P)}$$ is a
semidistance on $J.$
\end{theorem}

{\it Proof:} We need to prove that: $d^s(A,B)+d^s(B,C) \geq d^s
(A,C).$ Then it is sufficient to show: $$\frac{\sup_{P\in K}
g_{AB}(P)}{\inf_{P\in K} g_{AB}(P)} \cdot \frac{\sup_{Q\in K}
g_{BC}(Q)}{\inf_{Q\in K} g_{BC}(Q)} \geq \frac{\sup_{R\in K}
g_{AC}(R)}{\inf_{R\in K} g_{AC}(R)}.$$ Denote by $\alpha$ the left
hand side term in the inequality above and remark that $$\alpha
\geq \frac{g_{AB}(P)}{g_{AB}(Q)} \cdot \frac{g_{BC}(Q)}{g_{BC}(P)}
= \frac{\frac{f(P,A)}{f(P,B}}{\frac{f(Q,A)}{f(Q,B)}}\cdot
\frac{\frac{f(P,B)}{f(P,C}}{\frac{f(Q,B)}{f(Q,C)}}=\frac{g_{AC}(P)}{g_{AC}(Q)},
\ \ \ \ \ \forall P,Q \in K.$$ This means $\alpha \cdot g_{AC}(Q)
\geq g_{AC}(P),$ for all $P,Q \in K.$ Therefore, $$\alpha \cdot
g_{AC}(Q) \geq \sup_{P \in K} g_{AC}(P), \ \ \ \ \ \forall Q \in
K,$$ which yields $$\alpha \cdot \inf_{Q \in K} g_{AC}(Q) \geq
\sup_{P \in K} g_{AC}(P).$$ We obtain $$\alpha \geq \frac{\sup_{R
\in K} g_{AC}(R)}{\inf_{R \in K} g_{AC}(R)}.$$ \qed

The authors address their thanks to Professors M. Anastasiei, D.
E. Blair,  O. Dragi\v{c}evi\'{c} and John D. McCarthy for many
useful discussions on Barbilian's metrization procedure. The
authors thank also to the editor and the referees for their many
useful remarks that improved the content and the presentation of
this paper.

{\small

\bibliographystyle{amsplain}

}

Authors' addresses:

Wladimir G. Boskoff -
Department of Mathematics and Computer Science,
Ovidius University, Constan\c{t}a,
Romania; boskoff@univ-ovidius.ro

Marian G. Ciuc\u{a} -
Department of Mathematics and Computer Science,
Ovidius University, Constan\c{t}a,
Romania; mgciuca@univ-ovidius.ro

Bogdan D.
Suceav\u{a} Department of Mathematics, California State
University,Fullerton, CA 92834-6850,
U.S.A.; bsuceava@fullerton.edu

\end{document}